\documentclass[11pt,dvips]{article}
\input amssym.def
\input amssym.tex
\usepackage{psfig}
\usepackage{epsfig}
\def\binom#1#2{{#1}\choose{#2}}

\def\slfrac#1#2{\hbox{\kern.1em %
 \raise.5ex\hbox{\the\scriptfont0 #1}\kern-.11em %
 /\kern-.15em\lower.25ex\hbox{\the\scriptfont0 #2}}}

\newcommand{\eqn}[1]{(\ref{#1})}
\newcommand{\hsp}{\hspace*{\parindent}}

\newcommand{\eeq}{\end{equation}}
\newcommand{\beql}[1]{\begin{equation}\label{#1}}
\newcommand{\bsq}{{\vrule height .9ex width .8ex depth -.1ex }}

\newcommand{\sH}{{\cal H}}

\newcommand{\ZZ}{{\Bbb Z}}
\newcommand{\RR}{{\Bbb R}}
\newcommand{\FF}{{\Bbb F}}

\newcommand{\CC}{{\Bbb C}}

\newcommand{\dd}{\ldots}

\newcommand{\sG}{{\cal G}}

\newcommand{\sC}{{\cal C}}

\newcommand{\sQ}{{\cal Q}}
\newcommand{\sT}{{\cal T}}

\makeatletter
\def\@sect#1#2#3#4#5#6[#7]#8{\ifnum #2>\c@secnumdepth
     \def\@svsec{}\else
     \refstepcounter{#1}\edef\@svsec{\csname the#1\endcsname.\hskip .75em }\fi
     \@tempskipa #5\relax
      \ifdim \@tempskipa>\z@
        \begingroup #6\relax
          \@hangfrom{\hskip #3\relax\@svsec}{\interlinepenalty \@M #8\par}%
        \endgroup
       \csname #1mark\endcsname{#7}\addcontentsline
         {toc}{#1}{\ifnum #2>\c@secnumdepth \else
                      \protect\numberline{\csname the#1\endcsname}\fi
                    #7}\else
        \def\@svsechd{#6\hskip #3\@svsec #8\csname #1mark\endcsname
                      {#7}\addcontentsline
                           {toc}{#1}{\ifnum #2>\c@secnumdepth \else
                             \protect\numberline{\csname the#1\endcsname}\fi
                       #7}}\fi
     \@xsect{#5}}
\def\@begintheorem#1#2{\it \trivlist \item[\hskip \labelsep{\bf #1\ #2.}]}

\def\plain{plain}\ifx\fmtname\plain\csname fi\endcsname
     
     \input docstrip
     \preamble

     Do not distribute the stripped version of this file.
     The checksum in the header refers to the documented version.

     \endpreamble
     \generateFile{here.sty}{t}{\from{here.doc}{}}
     \endinput
\fi
\ifcat a\noexpand @\let\next\relax\else\def\next{%
    \documentstyle[here,doc]{article}\MakePercentIgnore}\fi\next
\ifx\@Hxfloat\@Hundef\else\expandafter\endinput\fi
\let\@Hxfloat\@xfloat
\def\@xfloat#1[{\@ifnextchar{H}{\@HHfloat{#1}[}{\@Hxfloat{#1}[}}
\def\@HHfloat#1[H]{%
\expandafter\let\csname end#1\endcsname\end@Hfloat
\vskip\intextsep\vbox\bgroup\def\@captype{#1}\parindent\z@
\ignorespaces}
\def\end@Hfloat{\egroup\vskip \intextsep}

\makeatother

\begin{document}
\begin{center}
{\Large {\bf A Family of Optimal Packings in Grassmannian Manifolds}} \\
\vspace{1.5\baselineskip}
{\em P. W. Shor} and {\em N. J. A. Sloane} \\
\vspace*{+.2\baselineskip}
Information Sciences Research \\
AT\&T Research\\
Murray Hill, NJ 07974 \\
\vspace{1\baselineskip}
Revised December 13, 1996 \\
\vspace{1.5\baselineskip}
{\bf ABSTRACT}
\vspace{.5\baselineskip}
\end{center}
\setlength{\baselineskip}{1.5\baselineskip}

A remarkable coincidence has led to the discovery of a family of
packings of $m^2 +m-2$
\linebreak
$m/2$-dimensional subspaces of $m$-dimensional
space, whenever $m$ is a power of 2.
These packings meet the ``orthoplex bound'' and are therefore optimal.

\vspace*{+3.5in}
\noindent\rule{1.25in}{.01in} 
\linebreak
\hspace*{+.85in}quantum coding theory, Clifford groups
\clearpage
\thispagestyle{empty}
\setcounter{page}{1}

\section{Introduction}
\hsp
Let $G(m,n)$ denote the Grassmannian space of all $n$-dimensional subspaces of real Euclidean $m$-dimensional space $\RR^m$.
The principal angles $\theta_1, \dd, \theta_n \in [0, \pi/2 ]$ between two
subspaces $P$, $Q \in G(m,n)$ are defined by
$$\cos \theta_i = \max_{u \in P} \max_{v \in Q} u \cdot v = u_i \cdot v_i ~,$$
for $i=1, \ldots, n$, subject to $u \cdot u = v \cdot v =1$,
$u \cdot u_j =0$, $v \cdot v_j =0$ $(1 \le j \le i-1)$.
We define the distance\footnote{It is shown in \cite{Grass} that this
is a metric, and in fact is essentially the $L_2$ distance
between the matrices that describe the orthogonal
projections onto $P$ and $Q$.}
between $P$ and $Q$ to be
$$d(P,Q) = \sqrt{\sin^2 \theta_1 + \cdots + \sin^2 \theta_n} ~.$$
In \cite{Grass} we discussed the problem of finding good packings in $G(m,n)$, that is, for given $N=1,2, \ldots$, of choosing
$P_1, \ldots, P_N \in G(m,n)$ so that $\min\limits_{i \neq j} d(P_i, P_j )$ is maximized.
It was shown that
for $N > m(m+1)/2$ the highest achievable distance, $d_N (m,n)$,
satisfies
\beql{Eq1}
d_N^2 (m,n) \le \frac{n(m-n)}{m} ~.
\eeq
A necessary condition for equality to hold in \eqn{Eq1} is that $N \le (m-1) (m+2)$.
An especially interesting case occurs when $m$ is even, $n=m/2$, and $N= (m-1) (m+2)$, where we found optimal packings for $m=2,4$ and 8;
that is,
packings of 4 lines in $\RR^2$,
18 2-spaces in $\RR^4$ and 70 4-spaces in $\RR^8$.
The first is the familiar configuration seen on the British flag (the Union Jack),
the second is the ``double-nine'', a classic configuration from nineteenth-century geometry (see the references in \cite{Grass}
and also \eqn{Eq4} below), but the third was discovered only after a very considerable computer-assisted search.
At the time \cite{Grass} was written we believed that there would be no further
examples in this series.

It came as a considerable surprise therefore when we discovered that such packings exist whenever $m$ is a power of 2.

These packings were discovered by a remarkable coincidence.
One of us (P.W.S.) had discovered a family
of groups in connection with quantum coding theory
\cite{CS96}, and asked the other (N.J.A.S.)
for the best way to determine their orders.
N.J.A.S. explained to P.W.S. that the Magma computer system
\cite{Mag1}, \cite{Mag2}, \cite{Mag3}
was ideal for this, and
gave as an example the symmetry group of above-mentioned set of 70
4-spaces in $\RR^8$,
an 8-dimensional group of order $2^7 8! = 5160960$. To our astonishment, the
first of his groups that P.W.S. tested turned out to be
(almost) exactly the same group.

The version of the group that arises from quantum coding in fact has the coordinates
in a slightly nicer order,
and produces the 70 planes as the orbit of the plane spanned by the first
four coordinate vectors.
With the help of our colleague R.~H. Hardin we verified that the next
three groups in the series produced packings meeting the bound in 16, 32 and 64 dimensions.
Further investigation then produced the general construction given in Section~3.
The groups are described in Section~2.

\section{The group}
\hsp
The group $\sG_i$ that arises from quantum coding theory
is a subgroup of the real orthogonal group $O(V, \RR )$, where $V$ denotes
$\RR^m$, $m=2^i$, $i \ge 1$, with coordinates
indexed by binary $i$-tuples $x= (x_1, \dd , x_i) \in \FF^{\,i}$, and $\FF$
is the field of order 2.
$\sG_i$ is generated by the following $2^i \times 2^i$ orthogonal matrices:

(i)~all permutation matrices $\pi_{A,b}$ corresponding to affine transformations $x \mapsto Ax + b$ of $\FF^{\,i}$, where $A$ is any invertible $i \times i$ matrix over $\FF$ and $b \in \FF^{\,i}$, and

(ii)~the matrix $H= {\rm diag} \{ H_2, H_2, \dd, H_2 \}$,
where $H_2 = \frac{1}{\sqrt{2}} {\binom{++}{+-}}$
(and $+$ denotes $+1$, $-$ denotes $-1$).

\noindent
By multiplying these generators it is easy to see that,
for $i \ge 2$,
$\sG_i$ contains the matrix $H' = {\rm diag} \{ H_4, H_4, \dd, H_4 \}$, where
$$H_4 = \frac{1}{2} \left( \matrix{
+ & + & + & + \cr
+ & - & + & - \cr
+ & + & - & - \cr
+ & - & - & + \cr
}
\right) ~.
$$
Let $\sH_i$ be the group generated by the permutations
and $H'$.
Then $\sG_i = \sH_i \bigcup H \sH_i$.

The packings described in Section~3 can be obtained by writing
the coordinates in the natural lexicographic order and taking the orbit
under $\sG_i$
of the subspace spanned by the first
$2^{i-1}$ coordinate vectors
(i.e., those in which $x_1 =0$).
However, the construction now given in Section~3 is a recursive one that no longer explicitly mentions the group.
The group is only needed in the analysis, where we make use of the fact that it acts transitively on the subspaces.
In the rest of this section we shall therefore give only a brief
discussion of these groups, in order to show their connection
with the Barnes-Wall lattices.

It turns out that $\sH_i$ and $\sG_i$
are well-known groups.
$\sH_i$ is the Clifford group $\sC \sT_1^+ (2^i )$ studied in
\cite{BRWI}, \cite{BRWII}, \cite{Wall}, which in recent years
has been used in the classification of finite simple groups (see the references
in \cite{CCKS}).
$\sH_i$ is relevant for the present work because of its connection
with the Barnes-Wall lattices.

Although the original paper of Barnes and Wall \cite{BW59} describes a family of lattices in each dimension $m=2^i$ $(i \ge 1)$, the most interesting lattices are the pair with the highest
number of minimal vectors (this number is given by the formula displayed in \eqn{Eq2}).
We denote this pair of $2^i$-dimensional lattices by $BW_i$ and $BW'_i$.
A construction of these lattices using
Reed-Muller codes is given in \cite{BS83} and in
\cite{SPLAG}, p.~234, example (f) (see also \cite{For88}).

$BW_i$ and $BW'_i$ are geometrically similar lattices, differing only by a rotation and change of scale.
When $i=1$, for example, we can take $BW_1$ to be the square lattice $\ZZ^2$
(Fig.~1, solid circles), and $BW'_1$ to be its sublattice of index 2
(Fig.~1, double circles).
In this case the matrix $D= \sqrt{2} H_2$ acts as an endomorphism sending $BW_1$ to $BW'_1$.
In exactly the same way, the matrix $\sqrt{2} H$ sends $BW_i$ to $BW'_i$,
a geometrically similar sublattice of
index $2^{m/2}$ (cf. \cite{SPLAG}, pp.~240--241).
Applying $\sqrt{2} H$ twice sends $BW_i$ to $2\,.\,BW_i$.
\begin{figure}[htb]
\centerline{\psfig{file=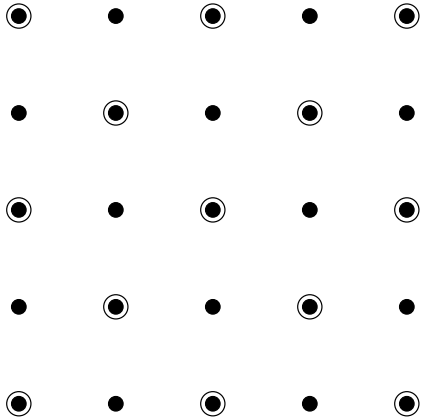,width=2.2in}}

\caption{The first pair of Barnes-Wall lattices, $BW_1$ (small circles) and $BW'_1$ (large circles).}
\end{figure}

Wall \cite{Wall} showed that for $i \neq 3$, $\sH_i$ is the full automorphism
group of the lattices $BW_i$ and $BW'_i$.
(The case $i=3$ is special, since $BW_3$ and $BW'_3$ are copies of the root lattice $E_8$.)
The group $\sH_i$ has a normal subgroup $E$ which is an extra-special 2-group
of order $2^{1+2i}$, and $\sH_i / E$ is isomorphic to the orthogonal
group $O_{2i}^+ (2) \cong D_i (2)$.
The order of $\sH_i$ is
\beql{Eq3}
2^{2i+1} \cdot 2^{i(i-1)} (2^i -1) \prod_{j=1}^{i-1} (4^j -1) ~.
\eeq
By adjoining the irrational matrix $H$ we obtain the full group
$\sG_i$, twice the size of $\sH_i$.
The group $\sG_i$ also appears in an apparently totally different context in \cite{CCKS} (see the group $L$ defined in Eq.~(2.13)).

The way the group $\sG_i$ arises in quantum coding theory is as follows.
The quantum state space of $i$ 2-state quantum systems is the complex space $\CC^m$, $m=2^i$.
Quantum computation involves making unitary transformations in this space (see \cite{CS96}, \cite{Shor95}).
Some transformations may be much easier to realize than others, and it is therefore important to know
which sets of transformations are sufficient for quantum computation, that is, which sets generate a group dense in $SU(2^i)$.
An interesting set of transformations which generate a finite group are the linear
Boolean functions on quantum bits (the permutation matrices in our group $\sG_i$),
and certain rotations of quantum bits by $\pi /2$.
To obtain the corresponding subgroup of the orthogonal group $SO(2^i)$, only one rotation is
required, which can be taken to be the matrix $H$.

\section{The construction}
\hsp
We specify a subspace $P\in G(m,n)$ by giving a generator matrix, that is,
an $n \times m$ matrix whose rows span $P$.
We will use the same symbol for the subspace and the generator matrix, and $P^\perp$ will denote the subspace orthogonal to $P$ (or a generator matrix
thereof).
$I$ denotes an identity matrix.

The construction is recursive.
We define a set $\sQ_i$ containing $2^{2i-1}$ monomial matrices of size
$2^{i-1} \times 2^{i-1}$ by $\sQ_1 = \{ (+), (-)\}$,
$$\sQ_i = \left\{
\left( \matrix{+ & 0 \cr 0 & + \cr} \right) \otimes Q ,~
\left( \matrix{+ & 0 \cr 0 & - \cr} \right) \otimes Q, ~
\left( \matrix{0 & + \cr + & 0 \cr} \right) \otimes Q , ~
\left( \matrix{0 & + \cr - & 0 \cr} \right) \otimes Q ; ~
Q \in \sQ_{i-1} \right\} ~,
$$
for $i \ge 2$.
Then $\sC_i$ is defined by
\begin{eqnarray*}
\sC_1 & = & \{ (+0) , (0+), (++), (+-) \} ~, \\
\sC_i & = & \left\{ (I0) , (0I) ,
\left( \matrix{P & 0 \cr 0 & P} \right) , \left(\matrix{P & 0 \cr 0 & P^\perp \cr} \right), (IQ) ;~
P \in \sC_{i-1}, ~Q \in \sQ_i \right\} ~,
\end{eqnarray*}
for $i \ge 2$.
For example, $\sC_2$ consists of the 18 matrices
\beql{Eq4}
\mbox{\footnotesize $\begin{array}{l}
\left( \begin{array}{c@{}c@{}c@{}c}
+& 0& 0& 0 \\ 0&+&0&0
\end{array} \right) ,~
\left(\begin{array}{c@{}c@{}c@{}c}
0&0&+&0 \\ 0&0&0&+
\end{array} \right) \\ [+.2in]
\left(\begin{array}{c@{}c@{}c@{}c}
+&0&0&0 \\ 0&0&+&0
\end{array}\right) ,~
\left(\begin{array}{c@{}c@{}c@{}c}
+&0&0&0 \\ 0&0&0&+
\end{array} \right),~
\left(\begin{array}{c@{}c@{}c@{}c}
0 &+ &0 &0 \\0& 0& 0& +
\end{array} \right),~
\left(\begin{array}{c@{}c@{}c@{}c}
0 &+ &0 &0\\0& 0& +& 0
\end{array} \right),~
\left(\begin{array}{c@{}c@{}c@{}c}
+ &+ &0 &0\\0& 0& +& +
\end{array} \right), \\ [+.2in]
\mbox{\hspace*{+1.58in}}\left(\begin{array}{c@{}c@{}c@{}c}
+ &+ &0 &0\\0& 0& +& - 
\end{array} \right),~
\left(\begin{array}{c@{}c@{}c@{}c}
+ &- &0 &0\\0& 0& +& -
\end{array} \right),~
\left(\begin{array}{c@{}c@{}c@{}c}
+ &- &0 &0\\0& 0& +& + 
\end{array} \right), \\ [+.2in]
\left(\begin{array}{c@{}c@{}c@{}c}
+ &0 &+ &0\\0& +& 0& + 
\end{array} \right),~
\left(\begin{array}{c@{}c@{}c@{}c}
+ &0 &+ &0\\0& +& 0& -
\end{array} \right),~
\left(\begin{array}{c@{}c@{}c@{}c}
+ &0 &- &0\\0& +& 0& +
\end{array} \right),~
\left(\begin{array}{c@{}c@{}c@{}c}
+ &0 &- &0\\0& +& 0& -
\end{array} \right),~
\left(\begin{array}{c@{}c@{}c@{}c}
+ &0 &0 &+\\0& +& +& 0
\end{array} \right),~ \\ [+.2in]
\mbox{\hspace*{+1.7in}}\left(\begin{array}{c@{}c@{}c@{}c}
+ &0 &0 &+\\0& +& -& 0
\end{array} \right),~
\left(\begin{array}{c@{}c@{}c@{}c}
+ &0 &0 &-\\0& +& +& 0
\end{array} \right),~
\left(\begin{array}{c@{}c@{}c@{}c}
+ &0 &0 &-\\ 0& +& -& 0 
\end{array} \right)\,.
\end{array}
$}
\eeq
(The third set of matrices in \eqn{Eq4} are the
matrices $(IQ)$.)
These are generator matrices for 18 2-spaces in $\RR^4$.
\paragraph{Theorem.}
{\em Let $m=2^i$, $i \ge 1$.
The generator matrices $\sC_i$ define a set of $(m-1)(m+2) =2^{2i} + 2^i -2$ $\slfrac{1}{2}m$-dimensional subspaces of $\RR^m$.
The distance between any two subspaces is
either
$\sqrt{m/4}$ or $\sqrt{m/2}$.
}
\paragraph{Proof.}
The number of subspaces is, by induction,
$$2+ 2(2^{2i-2} + 2^{i-1} -2) + 2^{2i-1} = 2^{2i} + 2^i -2 ~,$$
as claimed.

Since the recursive definition of the $\sC_i$ mentions the matrices $(I0)$ and $(0I)$, the coordinate
positions of $\sC_i$ can be labeled from left to right with binary $i$-tuples in the natural order, and the group
$\sG_i$ then acts by multiplication on the right.
It is now easy
to find matrices in $\sG_i$ that permute the subspaces transitively.
We leave the details to the reader.
Therefore, to determine the distances between the planes,
we may assume that one of the planes has generator matrix
$$A = \left[
\begin{array}{ccccclccccc}
1 & ~ & ~ & ~ & ~ && 0 & ~ & ~ & ~ & ~ \\
~ & 1 & ~ & ~ & ~ && ~ & 0 \\
~ & ~ & 1 & ~ & ~ && ~ & ~ & 0 \\
~ & ~ & ~ & \cdots &~ && ~ & ~ & ~ & \cdots \\
~ & ~ & ~ & ~ & 1 && ~ & ~ & ~ & ~ & 0
\end{array}
\right]
~.
$$
We recall (cf. \cite{Grass}) that if a second plane has generator matrix
$$B = \left[
\begin{array}{ccccclccccc}
c_1 & ~ & ~ & ~ & ~ & & s_1 \\
~ & c_2 & ~ & ~ & ~ & & ~ & s_2 \\
~ & ~ & c_3 & ~ & ~ & & ~ & ~ & s_3 \\
~ & ~ & ~ & \cdots & ~ & & ~ & ~ & ~ & \cdots \\
~ & ~ & ~ & ~ & c_n & & ~ & ~ & ~ & ~ & s_n
\end{array}
\right] ~,
$$
where $c_1^2 + s_1^2 = \cdots = c_n^2 + s_n^2 =1$, $n = 2^{i-1}$, then
the principal angles between $A$ and $B$
are $\arccos c_1$, $\arccos c_2, \dd, \arccos c_n$.

The principal angles between $A$ and $(0I)$
are $\pi /2$ ($n$ times).
Between $A$ and the subspaces
$$
\left( \matrix{ P & 0 \cr 0 & P \cr} \right) ~~~\mbox{or}~~~
\left( \matrix{ P & 0 \cr 0 & P^\perp \cr} \right)
$$
they are 0 ($n/2$ times), $\pi /2$ ($n/2$ times);
and between $A$ and $(IQ)$ they are $\pi /4$ ($n$ times).
The distance from $A$ to any other plane is therefore either
$\sqrt{n/2}$ or $\sqrt{n}$. \hfill $\bsq$

Since the bound \eqn{Eq1} is achieved, this is an optimal packing.

Together with R.~H. Hardin, we are also investigating other families of
subspaces that can be obtained from the same group.
If the initial subspace is taken to be that spanned by the first coordinate vector,
the orbit consists of the minimal vectors of the Barnes-Wall lattice $BW_i$,
together with their images under the transformation $H$, giving a total of
\beql{Eq2}
(2+2) (2^2 +2) \cdots (2^i +2)
\eeq
lines, with minimal angle $\pi /4$.
Taking the plane spanned by the first two coordinates as the initial plane, we
appear to obtain packings in $G(m,2)$ containing
$$\frac{1}{12} (2^i -1) \prod_{r=0}^i (2^r +2)$$
planes, with minimal distance 1, for $m=2^i$, $i \ge 1$.

On the other hand, if the initial subspace is that generated by the
first $m/4$ coordinate vectors, we appear to obtain packings in
$G(m,m/4)$ containing
$$\frac{1}{12} (m-2) (m-1) (m+2) (m+4)$$
subspaces, with minimal distance $\sqrt{m/8}$,
for $m=2^i$, $i \ge 2$.
The first member of this sequence is the packing of 24 lines in $\RR^4$
formed from the diameters of a pair of dual 24-cells.

We hope to discuss these packings
(which appear to be a kind of Grassmannian analogue of Reed-Muller
codes and Barnes-Wall lattices) in a subsequent paper.

\subsection*{Acknowledgements}
\hsp
We thank J. H. Conway and W. M. Kantor for helpful discussions
about the groups and R.~H. Hardin for assistance with experiments on the computer.

\end{document}